\documentclass[12pt,leqno,a4paper]{article}
\usepackage{amsmath}
	\allowdisplaybreaks
\usepackage{amsthm,amssymb,mathrsfs}
\usepackage[dvipsnames]{xcolor}
\usepackage{stmaryrd}
\usepackage{dsfont}
\usepackage{mathtools}
	
	\newtagform{bold}[\textbf]{\bf (}{\bf )}
	\usetagform{bold}
\usepackage[shortlabels]{enumitem}
	\setlist[enumerate]{resume}
	\setlist[enumerate]{topsep=0em,itemsep=.5em,partopsep=0em,parsep=0em}	
\usepackage{graphics}
\usepackage[margin=3cm]{geometry}
\usepackage{changepage}
\usepackage{multicol}

\usepackage[colorlinks=true,
    bookmarksnumbered=true, linktocpage=true,
	 linkcolor = blue, citecolor = magenta, urlcolor = magenta]{hyperref}
\usepackage{bookmark}

\usepackage{scrextend}
\usepackage{extarrows}
\usepackage{tikz}

\usepackage{tikz-cd}
	\tikzcdset{every label/.append style = {font = \small}}	

\usepackage{titlesec}
	\titleformat{\section}{\bf\filcenter}{\thesection.}{.5em}{}
	\titleformat{\subsection}{\bf}{\thesubsection}{0em}{}
	\titlespacing*{\section}{0pt}{1em}{.5em}
	\titlespacing*{\subsection}{0pt}{.5em}{.25em}

\usepackage{tocloft}

	\cftsetindents{subsection}{1.5em}{1em}

\theoremstyle{definition}

\numberwithin{equation}{section}

\newtheorem{theorem}[equation]{Theorem}

\title{\large On integral images of Curtis homomorphisms for $\mathrm{GL}_n$ and $\mathrm{U}_n$}
\author{\normalsize Tzu-Jan Li\footnote{Institute of Mathematics, Academia Sinica, 6F, Astronomy-Mathematics Building, No.\;1, Sec.\;4, Roosevelt Road, Taipei 106319, TAIWAN. Email: {\tt tzujan@gate.sinica.edu.tw}}}
\date{}
\begin{document}

\maketitle

\vspace{-15mm}
\begin{abstract}
\noindent{\bf Abstract.} For $G=\mathrm{GL}_n$ or $\mathrm{U}_n$ defined over a finite field of characteristic $p$, we refine a result of Bonnaf{\'e} and Kessar on the saturatedness of the Curtis homomorphism $\mathrm{Cur}^G$ by describing the image of $\mathrm{Cur}^G$ over $\overline{\mathbb{Z}}[1/p]$ via a system of linear conditions.\\[-3mm]

\end{abstract}

\section{Introduction}\label{intro-section}

\indent\indent Let $\mathbb{F}_q$ be a finite field with $q$ elements and of characteristic $p$, let $G$ be a connected reductive group defined over $\mathbb{F}_q$, and let $F:G\longrightarrow G$ be the associated Frobenius endomorphism, so $G^F=\{g\in G=G(\overline{\mathbb{F}}_q):F(g)=g\}=G(\mathbb{F}_q)$ is a finite group.

Let $\Lambda$ be a subring of $\overline{\mathbb{Q}}$ containing $\overline{\mathbb{Z}}[1/p]$, and then consider $\mathbf{E}_G=\mathrm{End}_{\Lambda[G^F]}(\Gamma_G)$, the endomorphism algebra of a Gelfand--Graev representation $\Gamma_G$ of $G^F$ with coefficients in $\Lambda$. It is known that $\mathbf{E}_G=\Lambda\mathbf{E}_G$ is a commutative $\Lambda$-algebra which is independent of the choice of $\Gamma_G$ up to isomorphism (see \cite[Sec.\;1.2--1.3]{Endomorphism}).

Let $\mathscr{T}_G=\mathscr{T}_{G,F}$ be the set of all $F$-stable maximal tori of $G$. For every $S\in\mathscr{T}_G$, Curtis has constructed  a $\overline{\mathbb{Q}}$-algebra homomorphism
\[
\mathrm{Cur}_S^G:\overline{\mathbb{Q}}\mathbf{E}_G \longrightarrow \overline{\mathbb{Q}}[S^F]
\]
compatible with the irreducible characters of $\overline{\mathbb{Q}}\mathbf{E}_G$ (see \cite[Th.\;4.2]{Curtis}). The homomorphism $\mathrm{Cur}_S^G$ is defined over $\Lambda$, in the way that $\mathrm{Cur}_S^G(\Lambda\mathbf{E}_G)\subset\Lambda[S^F]$ (see \cite[Lem.\;1.5(a)]{Endomorphism}). We may then form the ``Curtis homomorphism''
\[
\mathrm{Cur}^G:=(\mathrm{Cur}_{S}^G)_{S\in\mathscr{T}_G}:\overline{\mathbb{Q}}\mathbf{E}_G\longrightarrow\prod_{S\in\mathscr{T}_G}\overline{\mathbb{Q}}[S^F],
\]
which is an injective $\overline{\mathbb{Q}}$-algebra homomorphism (see \cite[Cor.\;3.3]{Bonnafe-Kessar}). Observe that
\begin{equation}\label{inclusion}
\mathrm{Cur}^G(\Lambda\mathbf{E}_G)\subset\mathrm{Cur}^G(\overline{\mathbb{Q}}\mathbf{E}_G)\cap\prod_{S\in\mathscr{T}_G}\Lambda[S^F].\\[-2mm]
\end{equation}

Let $W$ be the Weyl group of $G$. Bonnaf{\'e} and Kessar have proved in \cite[Th.\;3.7]{Bonnafe-Kessar} that the inclusion (\ref{inclusion}) is an equality (in other words, $\mathrm{Cur}^G$ is ``saturated over $\Lambda$'') if $|W|^{-1}\in\Lambda$. When $|W|^{-1}\not\in\Lambda$, the inclusion (\ref{inclusion}) can be strict (see \cite[3.29]{Thesis} or \cite[Rmk.\;3.9]{Bonnafe-Kessar}), and in this case it is natural to ask how to describe the $\Lambda$-lattice $\mathrm{Cur}^G(\Lambda\mathbf{E}_G)$ in $\prod_{S\in\mathscr{T}_G}\Lambda[S^F]$. 

The main goal of this article is to give a description of $\mathrm{Cur}^G(\Lambda\mathbf{E}_G)$ when $G$ is a general linear group $\mathrm{GL}_n$ or a unitary group $\mathrm{U}_n$, without assuming the invertibility of $|W|$ in $\Lambda$. To do this, for every $F$-stable Levi subgroup $L$ of $G$, let us fix a choice of quasi-split $T_L\in\mathscr{T}_L$, set $W_L=N_L(T_L)/T_L$ to be the Weyl group of $L$, and choose for each $w\in W_L$ a corresponding $T_{L,w}\in\mathscr{T}_L$ of type $w$ relative to $T_L$, so $T_{L,w}=xT_Lx^{-1}$ for some $x\in L$ with $x^{-1}F(x)T_L=w$. Then our main result is:

\begin{theorem}\label{mainthm}
{\it
Let $G$ be $\mathrm{GL}_n$ or $\mathrm{U}_n$ (defined over $\mathbb{F}_q$) with $n\in\mathbb{Z}_{>0}$, and let $\Lambda$ be a subring of $\overline{\mathbb{Q}}$ containing $\overline{\mathbb{Z}}[1/p]$. Then
\[
\mathrm{Cur}^G(\Lambda\mathbf{E}_G)=\mathrm{Cur}^G(\overline{\mathbb{Q}}\mathbf{E}_G)\cap\Omega
\]
where $\Omega$ is the set of the elements $(f_S)_{S\in\mathscr{T}_G}$ of $\prod_{S\in\mathscr{T}_G}\Lambda[S^F]$ such that
\[
\frac{1}{|W_L|}\sum_{w\in W_L}(-1)^{\ell(w)}f_{T_{L,w}}(s)\in\Lambda
\]
for every $F$-stable Levi subgroup $L$ of $G$ and every $s\in Z(L)^F$. Here, $\ell(w)$ is the length of $w\in W_L$ (defined through the simple reflections of $W_L$), and $Z(L)$ is the centre of $L$.\\[-7mm]
}
\end{theorem}

Our Theorem \ref{mainthm}  refines \cite[Th.\;3.7]{Bonnafe-Kessar} for $G=\mathrm{GL}_n$ or $\mathrm{U}_n$. Indeed, if $|W|$ is invertible in $\Lambda$, then for all $F$-stable Levi subgroups $L$ of $G$ we have $|W_L|^{-1}\in\Lambda$  (since $|W_L|$ divides $|W|$), so $\Omega$ is the whole $\prod_{S\in\mathscr{T}_G}\Lambda[S^F]$ and hence Theorem \ref{mainthm} implies that (\ref{inclusion}) is an equality for $G=\mathrm{GL}_n$ or $\mathrm{U}_n$.

For $G=\mathrm{GL}_2$, Theorem \ref{mainthm} has been proved in \cite[Prop.\;3.27]{Thesis} by direct calculations. We remark that it is also possible to prove Theorem \ref{mainthm} for $G=\mathrm{GL}_3$ by similar (but much longer) calculations, while it seems to be difficult to proceed such direct calculations for $G=\mathrm{GL}_n$ with $n\geq 4$. 

In order to prove Theorem \ref{mainthm} in its full generality, the idea is to use a $\Lambda$-algebra isomorphism (a ``Fourier transform'')
\[
\Lambda\mathbf{E}_G\simeq\Lambda\mathbf{K}_{G^\ast}
\]
in \cite{Endomorphism2} (see also \cite[Th.\;10.1(1)]{Helm} and \cite[Th.\;3.13]{Endomorphism}) to show that Theorem \ref{mainthm} is a consequence of a theorem on the $\mathbf{K}_{G^\ast}$-side (Theorem \ref{K-thm0}), and then prove the latter theorem. Here, $G^\ast$ is a Deligne--Lusztig dual of $G$ (see \cite[Def.\;5.21]{Deligne--Lusztig}), and $\mathbf{K}_{G^\ast}$ denotes the Grothendieck ring of the category of $\overline{\mathbb{F}}_q[G^{\ast}(\mathbb{F}_q)]$-modules of finite $\overline{\mathbb{F}}_q$-dimension. Our argument will be based on the Jordan decomposition of characters, a tool which is well-adapted under our assumptions for $G$, but which will become delicate without these assumptions, mainly due to the existence of semisimple centralisers which are not Levi subgroups of $G$. The author hopes that Theorems \ref{mainthm} and \ref{K-thm0} can eventually be generalised to other reductive groups $G$, while new tools or new viewpoints may be needed for this purpose.

\section{A ``Fourier transform''}

\indent\indent In this section, $G$ is a connected reductive group defined over $\mathbb{F}_q$, and $F:G\longrightarrow G$ is the associated Frobenius endomorphism. Let $(G^\ast,F^\ast)$ be the dual of $(G,F)$ in the sense of Deligne and Lusztig. For a finite group $A$ and a field $k$, we shall denote by $\mathrm{Irr}_{k}(A)$ the set of irreducible characters of $A$ with values in $k$.

In \cite[Sec.\;2.5]{Endomorphism}, we have constructed a $\overline{\mathbb{Q}}$-algebra isomorphism  $\overline{\mathbb{Q}}\mathbf{E}_G\simeq\overline{\mathbb{Q}}\mathbf{K}_{G^\ast}$ satisfying the following property: For every $S\in\mathscr{T}_{G,F}$ and every $S^\ast\in\mathscr{T}_{G^\ast,F^\ast}$ dual to $S$, let $\mathrm{Res}_{S^{\ast F^\ast}}^{G^{\ast F^\ast}}:\mathbf{K}_{G^\ast}\longrightarrow\mathbf{K}_{S^\ast}$ be the restriction map and let $h:\mathbb{Z}[S^F]\xrightarrow{\;\sim\;}\mathbf{K}_{S^\ast}$ be the ring isomorphism induced by the toric duality $S^F\simeq\mathrm{Irr}_{\overline{\mathbb{Q}}}(S^{\ast F^\ast})$ (with respect to a fixed choice of identifications $(\mathbb{Q}/\mathbb{Z})_{p'}\simeq \overline{\mathbb{F}}_q^\times\hookrightarrow\overline{\mathbb{Q}}^\times$); then the following diagram of  $\overline{\mathbb{Q}}$-algebras is commutative:
\[
\begin{tikzcd}
\overline{\mathbb{Q}}\mathbf{E}_G\arrow[r,leftrightarrow,"\sim"]\arrow[d,"\mathrm{Cur}_S^G"'] & \overline{\mathbb{Q}}\mathbf{K}_{G^\ast}\arrow[d,"\mathrm{Res}_{S^{\ast F^\ast}}^{G^{\ast F^\ast}}"]\\
\overline{\mathbb{Q}}[S^F]\arrow[r,"h","\sim"'] & \overline{\mathbb{Q}}\mathbf{K}_{S^\ast}
\end{tikzcd}
\]
For every $f\in\overline{\mathbb{Q}}\mathbf{E}_G \simeq\overline{\mathbb{Q}}\mathbf{K}_{G^\ast}$ and every $s\in S^F$, a direct calculation shows that
\begin{equation}\label{dictionary}
\mathrm{Cur}_S^G(f)(s)=\langle f|_{S^{\ast F^\ast}},\widehat{s}\,\rangle_{S^{\ast F^\ast}},
\end{equation}
where $\langle\cdot,\cdot\rangle_{S^{\ast F^\ast}}$ is the standard pairing $\langle a,b\rangle_{S^{\ast F^\ast}}=|S^{\ast F^\ast}|^{-1}\sum_{s\in S^{\ast F^\ast}}a(s^{-1})b(s)$ for all $a,b\in\overline{\mathbb{Q}}\mathbf{K}_{S^\ast}$, and $\widehat{s}\in\mathrm{Irr}_{\overline{\mathbb{Q}}}(S^{\ast F^\ast})$ is the character corresponding to $s\in S^F$ by duality. On the other hand, for $\Lambda$ being a subring of $\overline{\mathbb{Q}}$ containing $\overline{\mathbb{Z}}[1/p]$ (as in Section \ref{intro-section}), from the study of \cite{Endomorphism2} we know that the $\overline{\mathbb{Q}}$-algebra isomorphism $\overline{\mathbb{Q}}\mathbf{E}_G\simeq\overline{\mathbb{Q}}\mathbf{K}_{G^\ast}$ here yields (by restriction) a $\Lambda$-algebra isomorphism
\begin{equation}\label{Fourier}
\Lambda\mathbf{E}_G\simeq\Lambda\mathbf{K}_{G^\ast}
\end{equation} 
whenever the bad prime numbers for $G$ (see {\it ibid.}) are all invertible in $\Lambda$. (If $G$ is as in Theorem \ref{mainthm}, there are no bad prime numbers and we can take $\Lambda=\overline{\mathbb{Z}}[1/p]$.) 

We would like to think of the isomorphism (\ref{Fourier}) as a ``Fourier transform,'' since it translates the convolution product of $\Lambda\mathbf{E}_G$ into the tensor product in $\Lambda\mathbf{K}_{G^\ast}$ (which corresponds to the pointwise product of Brauer characters).

Let $\mathbf{P}_{G^\ast}$ be the additive Grothendieck group of the category of projective $\overline{\mathbb{F}}_q[G^{\ast F^\ast}]$-modules of finite $\overline{\mathbb{F}}_q$-dimension, and view $\mathbf{P}_{G^\ast}$ as an ideal of $\mathbf{K}_{G^\ast}$. We then have a perfect pairing $\langle\cdot,\cdot\rangle$ (over $\mathbb{Z}$) between $\mathbf{K}_{G^\ast}$ and $\mathbf{P}_{G^\ast}$ defined by
\[
\langle V_1,V_2\rangle=\dim_{\overline{\mathbb{F}}_q}\mathrm{Hom}_{\overline{\mathbb{F}}_q[G^{\ast F^\ast}]}(V_1,V_2)
\]
where one of $V_1$ and $V_2$ belongs to $\mathbf{K}_{G^\ast}$ and the other belongs to $\mathbf{P}_{G^\ast}$ (see \cite[Sec.\;14.5]{Serre}). Moreover, denoting by $\mathrm{St}_{G^\ast}$ the Steinberg character of $G^{\ast F^\ast}$, the multiplication by (the reduction modulo $p$ of) $\mathrm{St}_{G^\ast}$ induces a $\mathbb{Z}$-module isomorphism from  $\mathbf{K}_{G^\ast}$ to $\mathbf{P}_{G^\ast}
$
(see \cite[Th.\;1.1]{Lusztig}). It follows that the following pairing is also perfect:
\begin{equation}\label{K-pairing}
\mathbf{K}_{G^\ast}\times\mathbf{K}_{G^\ast}\longrightarrow \mathbb{Z},\quad (V_1,V_2)\longmapsto \langle V_1,\mathrm{St}_{G^\ast}\cdot V_2\rangle.
\end{equation}

\section{Translation of Theorem \ref{mainthm} into the dual side}\label{K-section}

\indent\indent We define $\mathscr{L}_G=\mathscr{L}_{G,F}$ to be the set of $F$-stable Levi subgroups of $G$. For each $L\in\mathscr{L}_G$, let $R_{L}^{G}$ be the Lusztig induction (see {\normalfont\cite[Ch.\;9]{Digne-Michel}}). We also employ the notation $\epsilon_H=(-1)^{\mathrm{rank}_{\mathbb{F}_q}(H)}$ for every reductive group $H$ over $\mathbb{F}_q$.

 Let $f\in\overline{\mathbb{Q}}\mathbf{E}_G\simeq\overline{\mathbb{Q}}\mathbf{K}_{G^\ast}$ and let $L\in\mathscr{L}_{G}$. Let also $s\in Z(L)^F$, and denote by $\widehat{s}:L^{\ast F^\ast}\longrightarrow\overline{\mathbb{Q}}^\times$ the linear character dual to $s$, where $L^\ast\in\mathscr{L}_{G^\ast}=\mathscr{L}_{G^\ast,F^\ast}$ is a dual of $L$. Moreover, for each virtual $\overline{\mathbb{Q}}[G^{\ast F^\ast}]$-module $V$ of finite $\overline{\mathbb{Q}}$-dimension, let $\overline{V}\in\mathbf{K}_{G^\ast}$ be its reduction modulo $p$. Then, for each $S\in\mathscr{T}_{L}=\mathscr{T}_{L,F}$ (so $s\in S^F$) with $S^\ast\in\mathscr{T}_{L^\ast}=\mathscr{T}_{L^\ast,F^\ast}$ dual to $S$, we have
\[
\mathrm{Cur}_{S}^G(f)(s)=
\langle f,\overline{\mathrm{Ind}_{S^{\ast F^\ast}}^{G^{\ast F^\ast}}(\widehat{s}\,|_{S^{\ast F^\ast}})}\rangle
=\epsilon_{L^\ast}\epsilon_{S^\ast}
\langle f,\overline{\mathrm{St}_{G^\ast}\cdot R_{S^\ast}^{G^{\ast}}(\widehat{s}\,|_{S^{\ast F^\ast}})}\rangle
\]
by (\ref{dictionary}), the Frobenius reciprocity and \cite[Prop.\;7.3]{Deligne--Lusztig}. For $W_L$, $T_{L,w}$ and $\ell(w)$ as in Theorem \ref{mainthm}, we thus have\\[-3mm]
\begin{equation}\label{dic-eq2}
\begin{aligned}
\frac{1}{|W_L|}\sum_{w\in W_L}(-1)^{\ell(w)}\mathrm{Cur}_{T_{L,w}}^G(f)(s)
&=\frac{1}{|L^F|}\sum_{S\in\mathscr{T}_L}\epsilon_L\epsilon_S|S^F|\mathrm{Cur}_{S}^G(f)(s)\\
&=\frac{1}{|L^F|}\sum_{S\in\mathscr{T}_L}|S^F|\langle f,\overline{\mathrm{St}_{G^\ast}\cdot R_{S^{\ast}}^{G^{\ast}}(\widehat{s}\,|_{S^{\ast F^\ast}})}\rangle\\
&=\frac{1}{|L^{\ast F^\ast}|}\sum_{S^\ast\in\mathscr{T}_{L^\ast}} \langle f,\overline{\mathrm{St}_{G^\ast}\cdot R_{L^{\ast}}^{G^{\ast}}(\widehat{s}\cdot |S^{\ast F^\ast}|R_{S^{\ast}}^{L^{\ast}}(1))}\rangle\\
&=\langle f,\overline{\mathrm{St}_{G^\ast}\cdot R_{L^{\ast}}^{G^{\ast }}(\widehat{s}\,)}\rangle,
\end{aligned}
\end{equation}
where the first equality holds since the map $S\in\mathscr{T}_L\longmapsto \epsilon_L\epsilon_S|S^F|\mathrm{Cur}_S^G(f)(s)\in\overline{\mathbb{Q}}$ is invariant under the $L^F$-conjugation on $\mathscr{T}_L$ and since $\epsilon_L\epsilon_{T_{L,w}}=(-1)^{\ell(w)}$ for $w\in W_L$, and where the last equality follows from \cite[(7.14.1)]{Deligne--Lusztig}.

By (\ref{dic-eq2}), the perfect pairing (\ref{K-pairing}) and the $\Lambda$-algebra isomorphism (\ref{Fourier}), we see that Theorem \ref{mainthm} is a corollary of the following theorem:

\begin{theorem}\label{K-thm0}
{\it
Let $G$ be as in Theorem \ref{mainthm}. Then
\[
\mathbf{K}_{G^\ast}=\sum_{L^\ast\in\mathscr{L}_{G^\ast}}\mathbb{Z}\cdot \overline{R_{L^\ast}^{G^\ast}(X(L^{\ast F^\ast}))},
\]
where $X(L^{\ast F^\ast})$ is the abelian group of $\overline{\mathbb{Q}}$-valued linear characters of $L^{\ast F^\ast}$.
}\\[-7mm]
\end{theorem}

As every element of $\mathbf{K}_{G^\ast}$ is the reduction modulo $p$ of a virtual $\overline{\mathbb{Q}}[G^{\ast F^\ast}]$-module (see \cite[Th.\;33]{Serre}), we find that Theorem \ref{K-thm0} is a corollary of the following theorem:

\begin{theorem}\label{K-thm}
{\it
In the setup of Theorem \ref{K-thm0}, we have
\[
\mathrm{Irr}_{\overline{\mathbb{Q}}}(G^{\ast F^\ast})\subset\sum_{L^\ast\in\mathscr{L}_{G^\ast}}\mathbb{Z}\cdot R_{L^\ast}^{G^\ast}(X(L^{\ast F^\ast})).\\[-2mm]
\]
}
\end{theorem}

{\it Remark.} In the case of $G=\mathrm{PGL}_2$ over $\mathbb{F}_q$ with $q$ odd (and with split Frobenius $F$), we have $G^\ast=\mathrm{SL}_2$ over $\mathbb{F}_q$. Direct calculations similar to that made in \cite[Prop.\;3.27]{Thesis} show that Theorem \ref{mainthm} still holds for $G$ here, while Theorem \ref{K-thm} does not hold for $G^\ast$ here, since the two irreducible characters of $G^{\ast F^\ast}=\mathrm{SL}_2(\mathbb{F}_q)$ of degree $(q+1)/2$, called $\chi_{\alpha_0}^{\pm 1}$ in \cite[Table 12.1]{Digne-Michel}, are not $\mathbb{Z}$-linear combinations of characters induced (in the sense of Deligne--Lusztig) from $F^\ast$-stable maximal tori of $G^\ast$. Thus Theorem \ref{K-thm} is in general a result stronger than Theorem \ref{mainthm} (for the corresponding $G^\ast$ and $G$).

\section{Proof of Theorems \ref{mainthm}, \ref{K-thm0} and \ref{K-thm}}\label{proof-section}

\indent\indent From now on, let $G$ be $\mathrm{GL}_n$ or $\mathrm{U}_n$ (defined over $\mathbb{F}_q$). From Section \ref{K-section}, we have the implications 
\[
\mbox{Theorem \ref{mainthm}}\Longleftarrow\mbox{Theorem \ref{K-thm0}}\Longleftarrow\mbox{Theorem \ref{K-thm}},
\]
so it is sufficient to prove Theorem \ref{K-thm}.

Let us prove Theorem \ref{K-thm}; our proof will be a modification of the proof of \cite[Th.\;11.7.3]{Digne-Michel}. As Theorem \ref{K-thm} is stated on the dual side $(G^\ast,F^\ast)$, we shall swap $(G,F)$ and $(G^\ast,F^\ast)$ to simplify the notation, so that we now need to prove
\[
\mathrm{Irr}_{\overline{\mathbb{Q}}}(G^F)\subset\sum_{L\in\mathscr{L}_{G}}\mathbb{Z}\cdot R_{L}^{G}(X(L^{F}))
\]
when $G$ is $\mathrm{GL}_n$ or $\mathrm{U}_n$ (note that $\mathrm{GL}_n$ and $\mathrm{U}_n$ are both self-dual). 

Let $\varphi\in\mathrm{Irr}_{\overline{\mathbb{Q}}}(G^{F})$. Then we can find a semisimple element $s$ of $G^{\ast F^\ast}$ such that $\varphi$ belongs to the geometric Lusztig series $\mathcal{E}(G^{ F},(s))$ associated to $(G^{ F},(s))$ (see \cite[Prop.\;11.3.2]{Digne-Michel}). Now set $L^\ast=C_{G^\ast}(s)\in\mathscr{L}_{G^\ast}$, and choose a dual $L\in\mathscr{L}_{G}$ of $L^\ast$. 
By the Jordan decomposition of irreducible characters (see \cite[Th.\;11.4.3(ii) and Prop.\;11.4.8(ii)]{Digne-Michel}), there is a $\varphi'\in\mathcal{E}(L^{F},1)$ such that
\begin{equation}\label{eq4.1}
\varphi=\epsilon_{G}\epsilon_{L}R_{L}^{G}(\widehat{s}\cdot\varphi'),
\end{equation}
where $\widehat{s}:L^F\longrightarrow\overline{\mathbb{Q}}^\times$ is the linear character dual to $s$.

We next apply the theory of almost characters to analyse the structure of $\varphi'$. Let $W_{L}=N_{L}(T_{L})/T_{L}$ be the Weyl group of $L$ with respect to a quasi-split $T_{L}\in\mathscr{T}_{L}$, and let $\widetilde{W}_{L}=W_{L}\rtimes \langle F\rangle$ where $\langle F\rangle$ is the finite cyclic subgroup of the automorphism group of $W_L$ generated by $F:W_{L}\longrightarrow W_{L}$ (induced from $F:G\longrightarrow G$). Denoting by $\mathrm{Irr}_{\mathbb{Q}}(W_{L})^F$ the set of $F$-invariant elements of $\mathrm{Irr}_{\mathbb{Q}}(W_{L})$ (note that $\mathrm{Irr}_{\overline{\mathbb{Q}}}(W_{L})=\mathrm{Irr}_{\mathbb{Q}}(W_{L})$ by \cite[Cor.\;1.15]{Springer}), to every $\chi\in\mathrm{Irr}_{\mathbb{Q}}(W_{L})^F$  we may associate an ``almost character''
\[
R_\chi=\frac{1}{|W_{L}|}\sum_{w\in W_{L}}\widetilde{\chi}(wF)R_{T_{L,w}}^{L}(1),
\]
where each $T_{L,w}\in\mathscr{T}_L$ is of type $w$ relative to $T_L$ (see the paragraph just above Theorem \ref{mainthm}), and where $\widetilde{\chi}\in\mathrm{Irr}_{\mathbb{Q}}(\widetilde{W}_L)$ is a choice of extension of $\chi$ (compare \cite[Sec.\;11.6--11.7]{Digne-Michel}, \cite[Sec.\;2]{LuSr} and \cite[Sec.\;7.3]{Carter}). Following the proofs of \cite[Th.\;11.7.2--11.7.3]{Digne-Michel}, the description of unipotent characters of $\mathrm{GL}_n(\mathbb{F}_q)$ and $\mathrm{U}_n(\mathbb{F}_q)$ by almost characters can be extended to our case of $L^F$, in the way that we have the following two properties:
\begin{enumerate}[(i)]
\item $\mathcal{E}(L^F,1)=\{\delta_\chi R_\chi\,|\,\chi\in\mathrm{Irr}_{\mathbb{Q}}(W_{L})^F\}$, where each $\delta_\chi\in\{\pm 1\}$ depends on $\chi$.
\item Each $R_\chi$ (where $\chi\in\mathrm{Irr}_{\mathbb{Q}}(W_{L})^F$) is a $\mathbb{Z}$-linear combination of the (virtual) characters $R_M^L(1)$ where $M\in\mathscr{L}_L$.
\end{enumerate}
By (i) and (ii), our $\varphi'$ may thus be expressed as
\begin{equation}\label{eq4.3}
\varphi'=\sum_{M\in\mathscr{L}_{L}}c_{M}R_{M}^{L}(1)\quad\mbox{for some }c_{M}\in\mathbb{Z}. \\[-2mm]
\end{equation}

We finally deduce from (\ref{eq4.1}) and (\ref{eq4.3}) that
\[
\varphi=\epsilon_{G}\epsilon_{L}\sum_{M\in\mathscr{L}_{L}}c_{M} R_{M}^{G}(\widehat{s\,}|_{M^F}),
\]
and this completes the proof of Theorem \ref{K-thm} since $\mathscr{L}_L\subset\mathscr{L}_G$. \qed

\section*{Acknowledgements} 
\indent\indent I thank Professor Jean-Fran{\c c}ois Dat and Cheng-Chiang Tsai for having helpful discussions with me. I also thank the Institute of Mathematics, Academia Sinica for providing me a postdoc position to realise this work.


\begin{thebibliography}{BoKe}
\bibitem[BoKe]{Bonnafe-Kessar} C.\;Bonnaf{\'e} and R.\;Kessar, {\it On the endomorphism algebras of modular Gelfand--Graev representations}, J.\;Alg.\;320 (2008)
\bibitem[Ca]{Carter} R.\;W.\;Carter, {\it On the Representation Theory of the Finite Groups of Lie Type over an
Algebraically Closed Field of Characteristic 0}, Encyclopaedia of Mathematical Sciences Vol.\;77, Algebra IX, Springer
(1996)
\bibitem[Cu]{Curtis} C.\;W.\;Curtis, {\it On the Gelfand--Graev Representations of a Reductive Group over a Finite Field}, J.\;Alg.\;157 (1993)
\bibitem[DeLu]{Deligne--Lusztig} P.\;Deligne and G.\;Lusztig, {\it Representations of reductive groups over finite fields}, Ann.\;Math.\;103-1 (1976)
\bibitem[DiMi]{Digne-Michel} 
F.\;Digne and J.\;Michel, Representations of finite groups of Lie Type, second edition, London Math.\;Soc.\;Student Texts 95, Cambridge Univ.\;Press (2020)
\bibitem[He]{Helm} D.\;Helm, {\it Curtis homomorphisms and the integral Bernstein center for $\mathrm{GL}_n$}, Alg.\;\& Num.\;Th.\;14-10 (2020)
\bibitem[Li1]{Thesis} T.-J.\;Li, {\it Sur l'alg{\`e}bre d'endomorphismes des repr{\'e}sentations de Gelfand--Graev
et le $\ell$-bloc unipotent de $\mathrm{GL}_2$ $p$-adique avec $\ell\neq p$}, Th{\`e}se de doctorat, Sorbonne Univ.\;(2022)
\bibitem[Li2]{Endomorphism} T.-J.\;Li, {\it On endomorphism algebras of Gelfand--Graev representations}, Represent.\;Theory 27 (2023)
\bibitem[LiSh]{Endomorphism2} T.-J.\;Li \& J.\;Shotton, {\it On endomorphism algebras of Gelfand--Graev representations II}, Bull.\;London Math.\;Soc.\;(2023)
\bibitem[Lu]{Lusztig}
G.\;Lusztig, {\it Divisibility of Projective Modules of Finite Chevalley Groups by the Steinberg Module}, Bull.\;London Math.\;Soc.\;8-2 (1976)
\bibitem[LuSr]{LuSr}
G.\;Lusztig and B.\;Srinivasan, {\it The Characters of the Finite Unitary Groups}, J.\;Alg.\;49 (1977)
\bibitem[Se]{Serre} J.-P.\;Serre, {\it Linear Representations of Finite Groups}, Graduate Texts in Mathematics 42, Springer (1977)
\bibitem[Sp]{Springer} T.\;A.\;Springer, {\it A Construction of Representations of Weyl Groups}, Inventiones Math.\;44
(1978)
\end{thebibliography}
\end{document}